\documentclass[12pt]{article}
\usepackage{amsmath, amsthm, amssymb}
\usepackage[dvips]{graphicx}
\usepackage{epsfig}

\newcommand{\A}{{\cal A}}

\newcommand{\eps}{{\varepsilon}}

\newcommand{\e}{\mathbb{E}}

\newcommand{\Reals}{\mathbb{R}}

\newcommand{\la}{\langle}
\newcommand{\ra}{\rangle}

\newtheorem{theorem}{Theorem}

\makeatletter
\@addtoreset{equation}{section}

\makeatother

\begin{document}

\title{The Parisi ultrametricity conjecture.}
\author{Dmitry Panchenko\thanks{Department of Mathematics, Texas A\&M University, email: panchenk@math.tamu.edu. 
Partially supported by NSF grant.}\\
}
\maketitle
\begin{abstract} 
In this paper we prove that the support of a random measure on the unit ball of a separable Hilbert space that
satisfies the Ghirlanda-Guerra identities must be ultrametric with probability one. This implies the Parisi ultrametricity 
conjecture in mean-field spin glass models, such as the Sherrington-Kirkpatrick and mixed $p$-spin models, 
for which Gibbs' measures are known to satisfy the Ghirlanda-Guerra identities in the thermodynamic limit.
\end{abstract}
\vspace{0.5cm}
\section{Introduction and main result.}

Let us consider a random probability measure $G$ on the unit ball of a separable Hilbert space $H$.
We will denote by $(\sigma^l)_{l\geq 1}$ an i.i.d. sample from this measure, by $\la\cdot\ra$
the average with respect to $G^{\otimes \infty}$ and by $\e$ the expectation with respect to the
randomness of $G$. Let $R_{l,l'} = \sigma^l\cdot \sigma^{l'}$ be the scalar product, or overlap,
of $\sigma^l$ and $\sigma^{l'}$. Random measure $G$ is said to satisfy the Ghirlanda-Guerra
identities if for any $n\geq 2,$ any bounded measurable function $f$ of the overlaps $(R_{l,l'})_{l,l'\leq n}$
and any bounded measurable function $\psi$ of one overlap,
\begin{equation}
\e \bigl\la f\psi(R_{1,n+1}) \bigr\ra = \frac{1}{n}\hspace{0.2mm} \e \bigl\la f \bigr\ra \hspace{0.2mm} \e \bigl\la \psi(R_{1,2}) \bigr\ra + \frac{1}{n}\sum_{l=2}^{n}\e \bigl\la f\psi(R_{1,l}) \bigr\ra.
\label{GG}
\end{equation}
Another way to express the Ghirlanda-Guerra identities is to say that, under the measure $\e G^{\otimes \infty}$,
conditionally on $R^n = (R_{l,l'})_{l,l' \leq n}$  the distribution of $R_{1,n+1}$ is given by  the mixture 
\begin{equation}
\frac{1}{n} \hspace{0.2mm}\mu + \frac{1}{n}\hspace{0.2mm} \sum_{l=2}^n \delta_{R_{1,l}},
\label{GGgen}
\end{equation}
where $\mu$ is the distribution of one overlap $R_{1,2}$ under $\e G^{\otimes 2}$. We will prove the following.
\begin{theorem}\label{ThUltra}
Under (\ref{GG}), the distribution of $(R_{l,l'})_{l,l' \geq 1}$ is ultrametric, i.e.
\begin{equation}
\e \bigl\la I\bigl(R_{1,2}\geq \min(R_{1,3},R_{2,3})\bigr) \bigr\ra=1.
\label{ultra}
\end{equation}
\end{theorem}\noindent
It is known (Theorem 2 in \cite{PGG}) that if $G$ satisfies the Ghirlanda-Guerra identities and if $q^*$ is the supremum of the support of $\mu$ then with probability one the support of $G$ belongs to the sphere of radius $\sqrt{q^*}$ in $H$. Therefore, (\ref{ultra}) means that with probability one over the choice of the random measure $G$, the distances in the Hilbert space $H$ between three independent replicas $\sigma^1,\sigma^2$ and $\sigma^3$ sampled from $G$ must satisfy the ultrametric inequality, 
\begin{equation}
\|\sigma^1-\sigma^2\|\leq \max \bigl(\|\sigma^1-\sigma^3\|,\|\sigma^2-\sigma^3\|\bigr).
\end{equation}
Examples of random measures satisfying the Ghirlanda-Guerra identities arise in several mean-field
spin glass models, such as the Sherrington-Kirkpatrick model \cite{SK} and mixed $p$-spin models, for which
measures $G$ are defined as the asymptotic analogues of the Gibbs measures in the thermodynamic limit
by way of the Dovbysh-Sudakov representation \cite{DS}.
Originally, the Ghirlanda-Guerra identities (\ref{GG}) were proved in \cite{GG} on average over the inverse temperature 
parameters and later, in a closely related formulation, by introducing a small perturbation term to the Hamiltonian  of the model (\cite{SG}, \cite{SG2}),
but in some cases can be proved in a strong sense without perturbation (\cite{PGGmixed}).  

The ultrametric structure of the overlap array $(R_{l,l'})$ appeared implicitly in the original work of G. Parisi in \cite{Parisi79}, \cite{Parisi}
in which the famous Parisi formula for the free energy in the Sherrington-Kirkpatrick model was discovered.
The fact that the particular form of the array $(R_{l,l'})$ suggested in \cite{Parisi79}, \cite{Parisi} encoded some definite physical 
properties of the Gibbs measure, including ultrametricity, was found during the subsequent interpretation of the Parisi solution 
in the work of M. M\'ezard, G. Parisi, N. Sourlas, G. Toulouse and M.A. Virasoro in \cite{M1}, \cite{M2} (see \cite{MPV} for more details). 
The Parisi formula for the free energy was proved rigorously in a celebrated work of M. Talagrand 
in \cite{TPF} following the breakthrough invention of the replica symmetry breaking interpolation scheme by F. Guerra in \cite{Guerra}, 
which gave a very strong indirect support to the entire Parisi ansatz including the ultrametricity conjecture.
More recently, several results providing some direct mathematical support to the ultrametricity conjecture were proved
under an additional technical assumption that the overlaps take only finitely many values, i.e. $R_{1,2} \in \{q_1,\ldots, q_k\}$  
with probability one for some non-random values $(q_l)_{l\leq k}$. The first such result was proved by L.-P. Arguin and
 M. Aizenman in \cite{AA} as a consequence of the Aizenman-Contucci stochastic stability property \cite{AC} of the Gibbs measures 
 in the mixed $p$-spin models. Inspired by \cite{AA}, the author proved 
 a similar result based on the Ghirlanda-Guerra identities in \cite{PGG} (see \cite{PGG2} for an elementary proof) 
and M. Talagrand gave a different proof in \cite{Tal-New}. Unfortunately, in the Sherrington-Kirkpatrick and 
mixed $p$-spin models one expects the distribution of the overlap to have a continuous component (\cite{MPV}), 
so the results in \cite{AA}, \cite{PGG} and \cite{Tal-New} were not directly applicable to these models. 

In this paper we deduce ultrametricity (\ref{ultra}) without any assumptions on the distribution of the overlap and,
as a result,  one can now give a more direct approach to the Parisi formula for the free energy 
in the mixed $p$-spin models (see \cite{PPF}). The proof of Theorem \ref{ThUltra} utilizes a new representation 
of the Ghirlanda-Guerra identities that appears in Theorem \ref{Th1} below, which can be viewed as a new invariance principle 
for random measures that satisfy (\ref{GG}). The idea behind this representation was originally motivated by the stability property 
proved in \cite{ACGG}, which unified the Aizenman-Contucci stochastic stability and the Ghirlanda-Guerra identities; however,  
the proof we give here is based only on the Ghirlanda-Guerra identities.

 \smallskip
  \textbf{Acknowledgement.} The author would like to thank Michel Talagrand for constant encouragement 
of the efforts that lead to this work.

\section{Invariance principles.}\label{Sec2}

In this section, we will first prove a new invariance property for random measures that satisfy the Ghirlanda-Guerra identities 
in Theorem {\ref{Th1}} and then deduce from it a modified  version of the invariance principle in Theorem \ref{Th2}
which will be used in the proof of Theorem \ref{ThUltra} in Section \ref{Sec3}. 
Given $n\geq 1$, consider $n$ bounded measurable functions 
$f_1,\ldots, f_n:  \Reals\to\Reals$
and define
\begin{equation}
F(\sigma,\sigma^1,\ldots,\sigma^n) = f_1(\sigma\cdot\sigma^1)+\ldots+f_n(\sigma\cdot\sigma^n).
\label{F1}
\end{equation}
For $1\leq l\leq n$ we define
\begin{equation}
F_l(\sigma,\sigma^1,\ldots,\sigma^n) = F(\sigma,\sigma^1,\ldots,\sigma^n)
 - f_l( \sigma\cdot\sigma^l)+ \e \la f_l(R_{1,2}) \ra
\label{F2}
\end{equation}
and for $l\geq n+1$ we define
\begin{equation}
F_l(\sigma,\sigma^1,\ldots,\sigma^n) =F(\sigma,\sigma^1,\ldots,\sigma^n).
\label{F3}
\end{equation}
The definition (\ref{F3}) for $l\geq n+1$ will not be used in the statement, but will appear in the proof of the next result.
Let us recall the notation $R^n = (R_{l,l'})_{l,l'\leq n}.$ 

\begin{theorem}\label{Th1}
Suppose (\ref{GG}) holds and let $\Phi$ be a bounded measurable function of $R^n.$ Then
\begin{equation}
\e \la\Phi \ra =
\e\Bigl\la
\frac{\Phi \exp \sum_{l=1}^{n} F_l(\sigma^l,\sigma^1,\ldots,\sigma^n)}
{\la\exp F(\sigma,\sigma^1,\ldots,\sigma^n)\ra_{\hspace{-0.3mm}\mathunderscore}^n}
\Bigr\ra,
\label{main}
\end{equation}
where the average $\la\cdot\ra_{\hspace{-0.3mm}\mathunderscore}$ in the denominator is in $\sigma$ only for fixed  $\sigma^1,\ldots, \sigma^n$
and the outside average of the ratio is in $\sigma^1,\ldots, \sigma^n$.
\end{theorem}\noindent
When $n=1,$ it is understood that $\Phi$ is a constant. Notice that one can easily recover the original Ghirlanda-Guerra 
identities from (\ref{main}) by taking $f_1 = t\psi$ and $f_2=\ldots = f_n=0$ and computing the derivative at $t=0.$

\smallskip\noindent\textbf{Proof.}
Without loss of generality, let us assume that $\Phi$ takes values in $[0,1]$ and suppose that
$|f_l|\leq L$ for $1\leq l\leq n$ for some large enough $L.$ For $t\geq 0$, let 
\begin{equation}
\varphi(t) = 
\e\Bigl\la
\frac{\Phi \exp \sum_{l=1}^{n} t F_l(\sigma^l,\sigma^1,\ldots,\sigma^n)}
{\la\exp t F(\sigma,\sigma^1,\ldots,\sigma^n)\ra_{\hspace{-0.3mm}\mathunderscore}^n}
\Bigr\ra.
\end{equation}
We will show that the Ghirlanda-Guerra identities (\ref{GG}) imply that this function is constant,
thus, proving the statement of the theorem, $\varphi(0)=\varphi(1).$ If for $k\geq 1$ we denote
$$
D_{n+k} = \sum_{l=1}^{n+k-1}F_l(\sigma^l,\sigma^1,\ldots,\sigma^n)
-(n+k-1) F_{n+k}(\sigma^{n+k},\sigma^1,\ldots,\sigma^n)
$$
then one can easily compute by induction that (recall (\ref{F3}) and that we average in $\sigma$ 
only in the denominator of (\ref{main}))
$$
\varphi^{(k)}(t) = 
\e\Bigl\la
\frac{\Phi D_{n+1}\ldots D_{n+k} 
\exp \sum_{l=1}^{n+k} t F_l(\sigma^l,\sigma^1,\ldots,\sigma^n)}
{\la\exp t F(\sigma,\sigma^1,\ldots,\sigma^n)\ra_{\hspace{-0.3mm}\mathunderscore}^{n+k}}
\Bigr\ra.
$$
First, let us notice that $\varphi^{(k)}(0)=0.$ Indeed, if we denote $\Phi' = \Phi D_{n+1}\ldots D_{n+k-1}$
then $\Phi'$ is a function of the overlaps $(R_{l,l'})_{l,l'\leq n+k-1}$ and 
\begin{align}
\varphi^{(k)}(0) 
& = \,
\e\Bigl\la \Phi' \Bigl(
\sum_{l=1}^{n+k-1}F_l(\sigma^l,\sigma^1,\ldots,\sigma^n)
-(n+k-1) F_{n+k}(\sigma^{n+k},\sigma^1,\ldots,\sigma^n)
\Bigr)
\Bigr\ra
\nonumber
\\
& = 
\sum_{j=1}^{n}
\e\Bigl\la \Phi' \Bigl(
\sum_{l\not = j,l=1}^{n+k-1}f_j(R_{j,l})
+ \e \la f_j(R_{1,2}) \ra
-(n+k-1) f_{j}(R_{j,n+k})
\Bigr)
\Bigr\ra
=0
\label{phizero}
\end{align}
by the Ghirlanda-Guerra identities (\ref{GG}) applied to each term $j$. Now, since $|F_l| \leq L n$
and $|D_{n+k}| \leq 2L (n+k-1)n$ we get
\begin{eqnarray*}
|\varphi^{(k)}(t)|
&\leq&
\Bigl(\prod_{l=1}^{k} 2L(n+l-1)n\Bigr) \,
\e\Bigl\la
\frac{\Phi 
\exp \sum_{l=1}^{n+k} t F_l(\sigma^l,\sigma^1,\ldots,\sigma^n)}
{\la\exp t F(\sigma,\sigma^1,\ldots,\sigma^n)\ra_{\hspace{-0.3mm}\mathunderscore}^{n+k}}
\Bigr\ra
\\
&=&
\Bigl(\prod_{l=1}^{k} 2L(n+l-1)n\Bigr) \,
\e\Bigl\la
\frac{\Phi 
\exp \sum_{l=1}^{n} t F_l(\sigma^l,\sigma^1,\ldots,\sigma^n)}
{\la\exp t F(\sigma,\sigma^1,\ldots,\sigma^n)\ra_{\hspace{-0.3mm}\mathunderscore}^{n}}
\Bigr\ra
\\
&=&
\prod_{l=1}^{k} (n+l-1)\,(2Ln)^k \, \varphi(t).
\end{eqnarray*}
Consider arbitrary $T>0.$ Again, using that $|F_l| \leq L n$ it is obvious that
$\varphi(t) \leq  e^{2L T n^2}$ for $0\leq t\leq T$ and, therefore,
$$
|\varphi^{(k)}(t)| \leq e^{2L T n^2} \frac{(n+k-1)!}{(n-1)!}\,  (2Ln)^k.
$$
By (\ref{phizero}) and Taylor's expansion
$$
|\varphi(t)-\varphi(0)| \leq \max_{0\leq s\leq t} \frac{|\varphi^{(k)}(s)|}{k!}t^k
\leq
e^{2L T n^2}  \frac{(n+k-1)! }{k! \,(n-1)!} (2Ln t)^k.
$$
Letting $k\to \infty$ we get that $\varphi(t)=\varphi(0)$ for $t<(2Ln)^{-1}.$
Therefore, for any $t_0<(2Ln)^{-1}$ we again have 
$\varphi^{(k)}(t_0)=0$ for all $k\geq 1$ and by Taylor's expansion for $t_0\leq t \leq T,$
$$
|\varphi(t)-\varphi(t_0)| \leq \max_{t_0\leq s\leq t} \frac{|\varphi^{(k)}(s)|}{k!}(t-t_0)^k
\leq
e^{2L T n^2}  \frac{(n+k-1)! }{k! \,(n-1)!} (2Ln (t-t_0))^k.
$$
Letting $k\to\infty$ proves that $\varphi(t) = \varphi(0)$ for $0\leq t< 2(2Ln)^{-1}.$
We can continue in the same fashion to prove this equality for all $0\leq t< T$ and 
note that $T$ was arbitrary.
\qed

\noindent
Let us write down a corollary of Theorem \ref{Th1} on which the proof of Theorem \ref{ThUltra} will be based.
Consider a finite index set $\A.$ Given $n\geq 1$ and configurations $\sigma^1,\ldots,\sigma^n,$ 
let $(B_\alpha)_{\alpha\in\A}$ be a partition of the Hilbert space $H$ such that for each $\alpha\in\A$ the indicator 
$I_{B_\alpha} = I(\sigma\in B_\alpha)$ is a measurable function of $R^n$ and $(\sigma\cdot\sigma^l)_{l\leq n}$ and let
\begin{equation}
W_\alpha=W_\alpha(\sigma^1,\ldots,\sigma^n)=G(B_\alpha).
\label{WA}
\end{equation}
Let us define a map $T$ by
\begin{equation}
W=(W_\alpha)_{\alpha\in\A}\to T(W) = 
\Bigl(\frac{\la I_{B_\alpha} \exp F(\sigma,\sigma^1,\ldots,\sigma^n )\ra_{\hspace{-0.3mm}\mathunderscore}}
{\la \exp F(\sigma,\sigma^1,\ldots,\sigma^n )\ra_{\hspace{-0.3mm}\mathunderscore}} \Bigr)_{\alpha\in\A}.
\label{TA}
\end{equation}
The following holds.
\begin{theorem}\label{Th2}
Under (\ref{GG}), for any bounded measurable function $\varphi:\Reals^{n^2}\times\Reals^{|\A|}\to \Reals$,
\begin{equation}
\e \bigl\la  \varphi(R^n, W) \bigr\ra
=
\e\Bigl\la
\frac{ \varphi(R^n,T(W)) \exp \sum_{l=1}^{n} F_l(\sigma^l,\sigma^1,\ldots,\sigma^n)}
{\la\exp F(\sigma,\sigma^1,\ldots,\sigma^n)\ra_{\hspace{-0.3mm}\mathunderscore}^n}
\Bigr\ra.
\label{nA}
\end{equation}
\end{theorem}\noindent
\textbf{Proof.} 
For each $\alpha\in\A$ let us take integer $n_\alpha\geq 0$ and let $m=n+\sum_{\alpha\in\A} n_\alpha.$ 
Let $(S_\alpha)_{\alpha\in\A}$ be any partition of $\{n+1,\ldots,m\}$ such that $|S_\alpha |=n_\alpha.$
Consider a continuous function $\Phi:\Reals^{n^2}\to\Reals$ and let
$\Phi' = \Phi(R^n) \prod_{\alpha\in\A}\varphi_\alpha$, where 
\begin{equation}
\varphi_\alpha = I \bigl(\sigma^l \in B_\alpha, \forall l\in S_\alpha \bigr).
\label{phialpha}
\end{equation}
Let $f_l$ for $l\leq n$ be as in (\ref{F1}) and $f_{n+1}=\ldots=f_m=0.$ Let us now apply Theorem \ref{Th1} 
with these choices of functions $\Phi'$ and $f_l$. First of all, integrating out the coordinates 
$(\sigma^l)_{l>n}$, the left hand side of (\ref{main}) can be written  as
\begin{equation}
\e \la \Phi'  \ra 
=
\e \Bigl\la \Phi(R^n) \prod_{\alpha\in\A} \varphi_\alpha \Bigr\ra
=
\e \Bigl\la\Phi(R^n) \prod_{\alpha\in\A} W_\alpha^{n_\alpha}(\sigma^1,\ldots,\sigma^n) \Bigr\ra,
\label{lhscor}
\end{equation}
where $W_\alpha$'s were defined in (\ref{WA}). Let us now compute the right hand side of (\ref{main}). 
Since $f_{n+1}=\ldots=f_m=0,$ the denominator will be
$\bigl\la\exp F(\sigma,\sigma^1,\ldots,\sigma^n)\bigr\ra_{\hspace{-0.3mm}\mathunderscore}^m$  and
\begin{equation}
\sum_{l=1}^{m} F_l(\sigma^l,\sigma^1,\ldots,\sigma^m) 
= 
\sum_{l=1}^{n} F_l(\sigma^l,\sigma^1,\ldots,\sigma^n)
+
\sum_{l=n+1}^{m} F(\sigma^l,\sigma^1,\ldots,\sigma^n).  
\label{numcor}
\end{equation}
Since the denominator does not depend on  $(\sigma^{l})_{l>n}$,
integrating in the coordinate $\sigma^l$ for $l\in S_\alpha$ will produce a factor 
$$
\la I_{B_\alpha} \exp F(\sigma,\sigma^1,\ldots,\sigma^n )  \ra_{\hspace{-0.3mm}\mathunderscore}.
$$
For each $\alpha\in \A$ we have $|S_\alpha| = n_\alpha$ such coordinates and, therefore, the right hand side of (\ref{main}) 
is equal to
\begin{equation}
\e\Bigl\la
\frac{\Phi(R^n) \exp \sum_{l=1}^{n} F_l(\sigma^l,\sigma^1,\ldots,\sigma^n)}
{\la\exp F(\sigma,\sigma^1,\ldots,\sigma^n)\ra_{\hspace{-0.3mm}\mathunderscore}^n}
\prod_{\alpha\in\A} \Bigl(
\frac{\la I_{B_\alpha} \exp F(\sigma,\sigma^1,\ldots,\sigma^n )\ra_{\hspace{-0.3mm}\mathunderscore}}
{\la \exp F(\sigma,\sigma^1,\ldots,\sigma^n )\ra_{\hspace{-0.3mm}\mathunderscore}}
\Bigr)^{n_\alpha}
\Bigr\ra.
\label{rhscor}
\end{equation}
Comparing with (\ref{lhscor}), recalling (\ref{TA}) and approximating a continuous function $\phi$ on $[0,1]^{|\A|}$ by polynomials 
we get (\ref{nA}) first for products $\Phi(R^n) \phi(W)$, then for continuous functions $\varphi(R^n,W)$ and then for arbitrary
bounded measurable functions.
\qed

\section{Proof of Theorem \ref{ThUltra}.} \label{Sec3}

We mentioned in the introduction that $G$ is concentrated on the sphere of radius $\sqrt{q^*}$ so all $\sigma$ below will be of length 
$\|\sigma\| = \sqrt{q^*}$. Consider a symmetric non-negative definite matrix $A=(a_{l,l'})_{l,l'\leq n}$ such that $a_{l,l} = q^*$ for $l\leq n$.
Given $\eps>0$, we will write $x\approx a$ to denote that $x\in (a-\eps,a+\eps)$ and $R^n\approx A$ to denote that 
$R_{l,l'} \approx a_{l,l'}$ for all $l\not = l' \leq n$ and, for simplicity of notation, we will keep the dependence of $\approx$ on $\eps$ implicit.
Below, the matrix $A$ will be used to describe a set of constraints such that the overlaps in $R^n$ can take values close to $A$, 
\begin{equation}
\e\bigl\la I(R^n \approx A)\bigr\ra >0,
\label{support}
\end{equation}
for a given $\eps>0$. Let us introduce the notation 
\begin{equation}
a_n^* = \max (a_{1,n},\ldots, a_{n-1,n}).
\end{equation}
The main step in the proof of Theorem \ref{ThUltra} is the following result which will be based on the invariance principle of Theorem \ref{Th2}.
\begin{theorem} \label{ThObs} Under (\ref{GG}), given $\eps>0$, if the matrix $A$ satisfies (\ref{support}) and $a_n^* +\eps < q^*$ then
\begin{equation}
\e\Bigl\la
I \Bigl(
R^n\approx A, R_{1,n+1}\approx a_{1,n},\ldots, R_{n-1,n+1}\approx a_{n-1,n}, R_{n,n+1} < a_n^* +\eps
\Bigr)
\Bigr\ra
>0.
\label{extend}
\end{equation}
\end{theorem}\noindent
Theorem \ref{ThObs} will be used in the following way. Suppose that the matrix $A$ is such that $a_n^*<q^*$ and $A$ is 
in the support of the distribution of $R^n$ under $\e G^{\otimes \infty}$ which means that (\ref{support}) holds for all $\eps>0$.
Since $a_n^*+\eps<q^*$ for small $\eps>0$,  (\ref{extend}) holds for all $\eps>0$.
Therefore, the support of the distribution of $R^{n+1}$ under $\e G^{\otimes \infty}$ intersects
the event in (\ref{extend}) for every $\eps>0$ and since the support is compact it contains a point in the set
\begin{equation}
\Bigl\{ A' : a_{l,l'}' = a_{l,l'} \mbox{ for } l, l' \leq n, a_{l,n+1}' = a_{l,n} \mbox{ for } l\leq n-1, a_{n,n+1}' \leq a_n^* \Bigr\}.
\label{Aplus}
\end{equation}

\smallskip \noindent
\textbf{Proof of Theorem \ref{ThObs}.} We will prove (\ref{extend}) by contradiction, so suppose that the left  hand side is equal to zero. 
We will apply Theorem \ref{Th2} with ${\cal A}=\{1,2\}$ and the partition 
$$
B_1 = \bigl\{\sigma : \sigma\cdot \sigma^n \geq a_n^* +\eps \bigr\},\, B_2 = B_1^c.
$$ 
Since we assume that $a_n^* +\eps < q^*$, the set $B_1$ contains a small neighborhood of $\sigma^n$ and on the event 
$\{R^n\approx A\}$ its complement $B_2=B_1^c$ contains small neighborhoods of $\sigma^1,\ldots, \sigma^{n-1}$ since 
$R_{l,n}<a_{l,n}+\eps \leq a_n^* +\eps$ and, thus, on this event for $\sigma^1,\ldots,\sigma^n$ in the support of $G$ 
the weights  $W_1=G(B_1), W_2 = G(B_2)=1-W_1$ are strictly positive. Then, (\ref{support}) implies that we can find 
$0<p<p'<1$ and small $\delta>0$ such that 
\begin{equation}
\e\bigl\la
I\bigl(
R^n\approx A, W_1 \in (p,p')
\bigr)
\bigr\ra \geq \delta.
\label{littlec}
\end{equation}
Let us apply Theorem \ref{Th2} and (\ref{nA}) with the above partition, the choice of 
\begin{equation}
\varphi(R^n,W) = I\bigl(R^n\approx A, W_1 \in (p,p')\bigr)
\end{equation}
and the choices of functions  $f_1=\ldots=f_{n-1}=0$ and $f_{n}(x) = t I(x\geq a_n^*+\eps)$ for $t\in\Reals$.  
The terms that appear on the right hand side of (\ref{nA}) will become 
\begin{eqnarray}
\sum_{l=1}^{n} F_l(\sigma^l,\sigma^1,\ldots,\sigma^n) 
& = &
\sum_{l=1}^{n-1} t  I(R_{l,n} \geq a_n^*+\eps)
+
t\hspace{0.3mm}  \e \bigl\la I(R_{1,2}\geq a_n^* + \eps)\bigr \ra
\nonumber
\\
& = &
t\hspace{0.3mm} \e \bigl\la I(R_{1,2}\geq a_n^* + \eps)\bigr \ra =: t\gamma
\nonumber
\end{eqnarray}
since, again, on the event $\{R^n\approx A\}$ the overlaps $R_{l,n}< a_{l,n}+\eps \leq a_n^*+\eps$ for $l\leq n-1$ and
\begin{equation}
\bigl \la \exp F(\sigma,\sigma^1,\ldots,\sigma^n) \bigr \ra_{\hspace{-0.3mm}\mathunderscore} 
=
\bigl \la \exp t I(\sigma\cdot \sigma^n \geq a_n^*+\eps)  \bigr\ra_{\hspace{-0.3mm}\mathunderscore} 
= \Delta_t(W)
\nonumber
\end{equation}
where  
\begin{equation}
\Delta_t(W) = G(B_1) e^t + G(B_2) =W_1 e^t + 1-W_1.
\label{DeltatW}
\end{equation}
If $W=(W_1,W_2),$ the map $T_t(W)$ corresponding to (\ref{TA}) can now be written as
\begin{equation}
T_t(W) = \Bigl(\frac{W_1 e^t}{\Delta_t(W) }, \frac{1-W_1}{\Delta_t(W)} \Bigr).
\label{TtW}
\end{equation}
Since $\Delta_t(W)\geq 1$ for $t\geq 0$, in this case equation (\ref{nA}) together with (\ref{littlec}) implies
\begin{eqnarray}
\delta &\leq &
\e\Bigl\la
\frac{I(R^n\approx A, (T_t(W))_1\in (p,p')) \,  e^{ t\gamma} }
{ \Delta_t(W)^n}
\Bigr\ra
\nonumber
\\
&\leq &
\e\Bigl\la
I\bigl(R^n\approx A, (T_t(W))_1\in (p,p') \bigr)  \, e^{ t\gamma} 
\Bigr\ra.
\label{littlec2}
\end{eqnarray}
In the average $\la\cdot\ra$ on the right hand side let us fix $\sigma^1,\ldots, \sigma^{n-1}$ and consider
the average with respect to $\sigma^n$ first. Clearly, on the event $\{R^n\approx A\}$ such average will be taken over the set
\begin{equation}
\Omega(\sigma^1,\ldots,\sigma^{n-1}) = \bigl\{\sigma: \sigma\cdot \sigma^l \approx a_{l,n} \mbox{ for } l\leq n-1\bigr\}.
\label{Omega}
\end{equation}
Let us look at the diameter of this set on the support of $G$. Suppose that with positive probability over the choice
of the measure $G$ and replicas  $\sigma^1,\ldots, \sigma^{n-1}$ from $G$ satisfying the constraints in $A$
(i.e. $R_{l,l'}\approx a_{l,l'}$ for $l,l'\leq n-1$)  we can find two points $\sigma',\sigma''$
in the support of $G$ that belong to the set $\Omega(\sigma^1,\ldots,\sigma^{n-1})$ and such that 
$\sigma'\cdot \sigma'' < a_n^* + \eps.$ This would then imply (\ref{extend})  since for 
$(\sigma^n,\sigma^{n+1})$ in a small neighborhood of $(\sigma',\sigma'')$ the vector 
$(\sigma^1,\ldots,\sigma^n,\sigma^{n+1})$ would belong to the event 
$$
\bigl\{
R^n\approx A, R_{1,n+1}\approx a_{1,n},\ldots, R_{n-1,n+1}\approx a_{n-1,n}, R_{n,n+1} < a_n^* +\eps
\bigr\}
$$
on the left hand side of (\ref{extend}). Since we assume that the left hand side of (\ref{extend}) is equal to zero,
we must have that with probability one  over the choice of the measure $G$ and replicas $\sigma^1,\ldots, \sigma^{n-1}$ 
satisfying the constraints in $A$ any two points $\sigma',\sigma''$ in the support of $G$ that belong 
to the set $\Omega(\sigma^1,\ldots,\sigma^{n-1})$ satisfy $\sigma'\cdot \sigma'' \geq a_n^* + \eps.$
Now, let us also recall that in (\ref{littlec2}) we are averaging over $\sigma^n$ that satisfy the condition $(T_t(W))_1 \in (p,p').$
If we fix any such $\sigma'$ in the support of $G$ that satisfies this condition and belongs to the set (\ref{Omega}), then the Gibbs
average in $\sigma^n$ will be taken over its neighborhood $B_1 = B_1(\sigma') = \{\sigma'': \sigma'\cdot \sigma'' \geq a_n^* +\eps\}$
of measure $W_1 = W_1(\sigma') = G(B_1(\sigma'))$ that satisfies $(T_t(W))_1\in (p,p').$
One can easily check that the map in (\ref{TtW}) satisfies $T_t^{-1} = T_{-t}$ and
using this for $(T_t(W))_1\in (p,p')$ implies that 
$$
W_1(\sigma') \in \Bigl\{\frac{q e^{-t}}{qe^{-t} + 1-q} : q\in (p, p') \Bigr\}
$$
and, thus, $W_1(\sigma') \leq (1-p')^{-1}e^{-t}$. This means that the average on the right hand side of (\ref{littlec2}) 
over $\sigma^n$ for fixed $\sigma^1,\ldots, \sigma^{n-1}$ is bounded by $(1-p')^{-1} e^{-t} e^{t\gamma }$ and, therefore, for $t\geq 0$
\begin{equation}
0<\delta  \leq
\e\bigl\la
I\bigl (R^n\approx A, (T_t(W))_1\in (p,p') \bigr)  e^{ t\gamma} 
\bigr\ra \leq (1-p')^{-1} e^{-t(1-\gamma)}.
\label{contra}
\end{equation}
Since $A$ satisfies (\ref{support}), $1-\gamma = \e\la I(R_{1,2}< a_n^* + \eps)\ra > 0$ and
letting $t\to+\infty$ in (\ref{contra}) we arrive at contradiction.
\qed

\medskip
\noindent \textbf{Proof of Theorem \ref{ThUltra}.} The proof is again by contradiction. 
Suppose that ultrametricity is violated in which case there exist $a<b\leq c <q^*$ such that the matrix
\begin{equation}
\left(
\begin{array}{ l c r }
  q^* & a & b \\
  a & q^* & c \\
  b &  c & q^* \\
\end{array}
\right)
\label{violated}
\end{equation}
is in the support of the distribution of $R^3$ under $\e G^{\otimes \infty}$ so it satisfies (\ref{support}) for every $\eps>0$.
In this case Theorem  \ref{ThObs} implies the following. 
Given any $n_1, n_2, n_3\geq 1$ and $n=n_1+n_2+n_3$ we can find a matrix $A$ in the support of the distribution of 
$R^n$ under $\e G^{\otimes \infty}$ such that  for some partition of indices $\{1,\ldots,n\} = I_1\cup I_2\cup I_3$ 
with $|I_j| = n_j$ we have $j\in I_j$ for $j\leq 3$ and 
\begin{enumerate}
\item[(a)] 
$a_{l,l'}\leq c$ for all $l\not = l'\leq n$,

\item[(b)] 
$a_{l,l'} = a$ if $l\in I_1, l'\in I_2$, $a_{l,l'} = b$ if $l\in I_1, l'\in I_3$ and $a_{l,l'} = c$ if $l\in I_2, l'\in I_3$.
\end{enumerate}
This can be proved by induction on $n_1,n_2,n_3.$ First of all, by the choice of the matrix (\ref{violated}) this holds for $n_1=n_2=n_3=1.$
Assuming the claim holds for some $n_1,n_2$ and $n_3$ with the matrix $A$, let us show how one can increase any of the $n_j$'s by one.
For example, let us assume for simplicity of notation that $n\in I_3$ and show that the claim holds with $n_3+1.$ 
Since $a_n^*\leq c<q^*$, we can use the comment below Theorem \ref{ThObs} to find a matrix $A'$ in the support of 
the distribution of $R^{n+1}$ under $\e G^{\otimes \infty}$ that belongs to the set (\ref{Aplus}).
Hence, $a_{l,l'}' \leq c$ for all $l\not = l' \leq n+1$ and $a_{l,n+1}' = a_{l,n}$ for $l\leq n-1$ so, in particular, $a_{l,n+1}' = b$ if $l\in I_1$ 
and  $a_{l,n+1}' = c$ if $l\in I_2$ which means that $A'$ satisfies the conditions (a), (b)  with $I_3$ replaced by $I_3\cup \{n+1\}.$
In a similar fashion, one can increase the cardinality of $I_1$ and $I_2$ which completes the induction.
Now, let $n_1=n_2=n_3=m$, find the matrix $A$ as above and find $\sigma^1,\ldots, \sigma^n$ on the sphere
of radius $\sqrt{q^*}$ such that $R_{l,l'} = a_{l,l'}$ for all $l,l'\leq n.$ Let $\bar{\sigma}^j$ be the barycenter of 
the set $\{\sigma^l : l\in I_j\}$. Condition (a) implies that  
$$
\|\bar{\sigma}^j\|^2 = \frac{1}{m^2}\sum_{l\in I_j} \|\sigma^l\|^2 + \frac{1}{m^2}\sum_{l\not = l'\in I_j} R_{l,l'}
\leq \frac{mq^* + m(m-1) c}{m^2}
$$
and condition (b) implies that $\bar{\sigma}^1\cdot \bar{\sigma}^2 = a$, $\bar{\sigma}^1\cdot \bar{\sigma}^3 = b$
and $\bar{\sigma}^2\cdot \bar{\sigma}^3 = c$. Therefore,
$$
\|\bar{\sigma}^2-\bar{\sigma}^3\|^2 = \|\bar{\sigma}^2\|^2 + \|\bar{\sigma}^3\|^2 
- 2\bar{\sigma}^2\cdot \bar{\sigma}^3 \leq \frac{2(q^* -c)}{m}
$$
and $0< b-a = \bar{\sigma}^1 \cdot \bar{\sigma}^3 -\bar{\sigma}^1\cdot \bar{\sigma}^2 \leq K m^{-1/2}$.
Letting $m\to\infty$, we arrive at contradiction.
\qed

\end{document}